\documentclass{elsart}

\usepackage{epsfig}

\usepackage[ansinew]{inputenc}
\usepackage{amsmath,amsfonts,latexsym,amssymb,amsbsy}
\usepackage{graphics}
\usepackage{graphicx}
\usepackage{epsfig}
\usepackage{float}

\usepackage{amsmath,amssymb}

\usepackage{pgfplots}

\usepackage{scalefnt}
\usepackage{color}

\usepackage[caption=false]{subfig}

\definecolor{marin}{rgb}   {0.,   0.3,   0.7}
\definecolor{rouge}{rgb}   {0.8,   0.,   0.}
\definecolor{sepia}{rgb}   {0.8,   0.5,   0.}
\usepackage[colorlinks,citecolor=marin,linkcolor=rouge,
            bookmarksopen,
            bookmarksnumbered
           ]{hyperref}

\usepackage{epsf}

\newcommand{\eps}{\varepsilon}

\newcommand{\R}{\mathbb{R}}

\newcommand{\T}{\mathbb{T}}

\numberwithin{equation}{section}

\newcommand{\QED}{\mbox{}\hfill \raisebox{-0.2pt}{\rule{5.6pt}{6pt}\rule{0pt}{0pt}}
          \medskip\par}

\begin{document}
\begin{frontmatter}

\journal{Applied Numerical Mathematics}

\title{High-order splitting methods for separable non-autonomous
parabolic equations}

\author{M. Seydao\u{g}lu$^{a,b}$} 
\ead{muasey@imm.upv.es} 
 and
\author{{S. Blanes$^{a}$}\corauthref{cor}}\corauth[cor]{Corresponding
author.} 
\ead{serblaza@imm.upv.es} $ \, $
\address{$^{a}$
Instituto de Matem\'{a}tica Multidisciplinar, Building 8G, second
floor,\\ Universitat Polit\`{e}cnica de
Val\`{e}ncia. 46022 Valencia, Spain.  \\
$^{b}$ Department of Mathematics, Faculty of Art and Science,
Mu\c{s} Alparslan University, 49100 Mu\c{s}, Turkey. }

\begin{abstract}

We consider the numerical integration of non-autonomous separable
parabolic equations using high order splitting methods with
complex coefficients (methods with real coefficients of order
greater than two necessarily have negative coefficients). We
propose to consider a class of methods
   that
allows us to evaluate all time-dependent
operators at real values of the time, leading to schemes which are
stable and simple to implement. If the system can be considered as
the perturbation of an exactly solvable problem  and the flow of
the dominant part is advanced using real coefficients, it is
possible to build highly efficient methods for these problems. We
show the performance of this class of methods on several numerical
examples and present some new improved schemes.

\end{abstract}


\end{frontmatter}

\section{Introduction}

We consider the numerical integration of non-autonomous separable
parabolic equations using high order splitting methods with
complex coefficients. This class of methods has been recently used
for the numerical integration of the autonomous case, showing good
performances \cite{blanes13oho,castella09smw,hansen09hos}.
Splitting methods with real coefficients of order greater than two
necessarily have negative coefficients and can not be used for
solving these problems \cite{blanes05otn,goldman96nos,sheng89slp,
suzuki91gto}. However, solutions with complex coefficients with
positive real part exist, and some of these methods can provide a
high performance in spite the equations have to be solved in the
complex domain.
Previous works with applications
among other in celestial mechanics and quantum mechanics where
splitting methods with complex coefficients are considered already
exist
\cite{bandrauk06cis,bandrauk91ies,chambers03siw,mclachlan02sm,prosen06hon,suzuki90fdo,suzuki91gto,suzuki95hep}.

A straightforward application of splitting methods with complex
coefficients to non-autonomous problems require the evaluation of
the time-dependent functions on the operators at complex times,
and the corresponding flows in the numerical scheme are, in
general, not well conditioned. In this work we propose to consider
a class of splitting methods in which one set of the coefficients belong to
the class of real and positive numbers. This can allow to evaluate
all time-dependent operators at real values of the time, leading
to schemes which are stable and simple to implement.

If the system can be considered as the perturbation of an exactly
solvable problem (or easy to numerically solve) and the flow of
the dominant part is advanced using the real coefficients, it is
possible to build highly efficient methods for these problems.

\subsection{The problem}


Let us consider the non-autonomous separable PDE
%
\begin{eqnarray} \label{eq:para}
\frac{d u}{dt}=  {   } A(t,u)  + {   } B(t,u), \qquad u(0)=u_0,
\end{eqnarray}
$u(x,t)\in\R^D$, and where the (possibly unbounded) operators
${   } A$, ${   } B$ and ${   } A+{   } B$ generate $C^0$ semi-groups
for positive $t$ over a finite or infinite Banach space. Equations
of this form are encountered in the context of {\em parabolic}
partial differential equations, an example being the inhomogeneous
non-autonomous {\em heat equation}
\begin{eqnarray} \label{eq.1}
\frac{\partial u}{\partial t}= \alpha(t)\Delta u + V(x,t) u, \quad
\mbox{or} \quad
 \frac{\partial u}{\partial t}= \nabla (a(x,t)\nabla u) + V(x,t) u
\end{eqnarray}
where $t\geq 0$, $x \in \R^d$ or $x \in \T^d$ and $\Delta$ denotes
the Laplacian with respect to the spatial coordinates, $x$.
Another example corresponds to reaction-diffusion equations of
the form
\begin{equation}\label{eq:NoAutNL}
 \frac{\partial u}{\partial t}  = D(t)\Delta u + {   } B(t,u),
\end{equation}
where $D(t)$ is a matrix of diffusion coefficients (typically a
diagonal matrix) and $B$ accounts for the reaction part. In
general, $A(t,u),B(t,u)$ can also depend on $x,\nabla,$ etc.,
which are omitted for clarity in the presentation.


For simplicity, we write the non-linear equation (\ref{eq:para}) in the
(apparently) linear form
\begin{eqnarray} \label{eq:para2}
 \frac{d u}{dt}=  L_{A(t,u)}u  + L_{B(t,u)}u,
\end{eqnarray}
where $L_{A},L_{B}$ are the Lie operators associated to $ {   } A, {   } B$,
i.e.
\[
    L_{A(t,u)} \equiv {   } A(t,u) \frac{\partial}{\partial u}, \qquad
  L_{B(t,u)} \equiv {   } B(t,u) \frac{\partial}{\partial u}
\]
which act on functions of $u$

If the problem is autonomous,
the formal solution is given by
$u(t)=\e^{t(L_{A(u)} + L_{B(u)})}u_0$,
which is a short way to write
\[
  u(t)=\e^{t(L_{A(u)} + L_{B(u)})}u_0 = \left.
  \sum_{k=0}^{\infty} \frac{t^k}{k!}\left({   } A(u) \frac{\partial}{\partial u} +
  {   } B(u) \frac{\partial}{\partial u} \right)^k u \
  \right|_{u=u_0}.
\]

If the subproblems 
\begin{eqnarray} \label{eq.Split_Aut}
 \frac{d u}{dt}=  A(u) \qquad \mbox{ and } \qquad
 \frac{d u}{dt}=  B(u)
\end{eqnarray}
 have exact solutions or can efficiently be numerically solved,
it is usual to consider splitting methods as
numerical integrators
. If we denote by $\e^{h
L_{A(u)}},\ \e^{h L_{B(u)}}$ the exact $h$-flows for each problem in
(\ref{eq.Split_Aut}) (and for a sufficiently small time step, $h$)
the simplest method within this class is the {\em Lie-Trotter
splitting}
\begin{equation}
\e^{h L_{A(u)}} \, \e^{h L_{B(u)}} \qquad \mbox{ or } \qquad \e^{h L_{B(u)}} \, \e^{h L_{A(u)}},
\end{equation}
which is a first order approximation in the time step to the
solution
, while the {\em symmetrized} version
\begin{equation} \label{eq:strang2}
 S(h)= \e^{h/2\, L_{A(u)}} \, \e^{h L_{B(u)}} \, \e^{h/2\, L_{A(u)}} \qquad \mbox{ or } \qquad
 S(h)= \e^{h/2\, L_{B(u)}} \, \e^{h L_{A(u)}} \, \e^{h/2\, L_{B(u)}}
\end{equation}
is referred to as {\em Strang splitting}, and is an approximation
of order $2$, i.e. $S(h)= \e^{h( L_{A(u) + B(u)})}+ \mathcal{O}(h^3)$. Upon
using an appropriate sequence of steps, high-order approximations
can be obtained as
\begin{eqnarray} \label{eq:splittingmethod}
\Psi(h) = \e^{h b_{m+1} L_B} \, \e^{h a_m L_A} \,\cdots \, \e^{h b_2
L_B} \, \e^{h a_1 L_A} \, \e^{h b_1 L_B},
\end{eqnarray}
%
and methods with real coefficients at any order can be obtained
\cite{creutz89hhm,suzuki90fdo,yoshida90coh}. However, as already
mentioned, splitting methods of order greater than two (with real
coefficients) have at least one of the coefficients $a_i$ negative
as well as at least one of the coefficients $b_i$ so,  the flows
$\e^{t L_A}$ and/or $\e^{t L_B}$ may not be well defined (this is
indeed the case, for instance, for the Laplacian operator) and
this prevents the use of methods which embed negative
coefficients. For this reason, exponential splitting methods of at
most order $p = 2$ have been considered up to recently.

%
%

In order to circumvent this order-barrier, the papers
\cite{castella09smw} and \cite{hansen09hos} simultaneously
presented a systematic analysis for a class of composition methods
with complex coefficients having positive real parts. Using this
extension from the real line to the complex plane, the authors of
\cite{castella09smw} and \cite{hansen09hos} built up methods of
orders $3$ to $14$ by considering a technique known as {\em
triple-jump composition}. More efficient high order methods are
obtained in \cite{blanes13oho}.

%

In this work we are interested, however, in the numerical
integration of the non-autonomous problem (\ref{eq:para}) where
the use of complex coefficients involve additional constraints as
we will see. A method of choice for solving numerically
(\ref{eq:para}) consists in advancing the solution alternatively
along  the exact (or numerical) solutions of the two problems
\begin{eqnarray} \label{eq.Split_tDep}
 \frac{d u}{dt}= A(t,u)  \qquad \mbox{ and } \qquad
 \frac{d u}{dt}= B(t,u) .
\end{eqnarray}
The exact flows are, in general, not known. This is the case, for
example, if
$[L_{A(t_i,u)},L_{A(t_j,u)}]=L_{A(t_i,u)}L_{A(t_j,u)} - L_{A(t_j,u)}L_{A(t_i,u)}\neq 0$
(and similarly for
$B(t)$). If the exact solution is not known, it can be replaced by
a sufficiently accurate numerical approximation.

This procedure is equivalent to take the time as two new
coordinates, $t_1,t_2$
\begin{equation}\label{}
  \left\{ \begin{array}{l}
  u'=  A(t_1, u) +  B(t_2, u)\\
 t_1' = 1\\
 t_2' = 1,
  \end{array} \right.
\end{equation}
with $'\equiv\frac{d}{dt}$, and to split the system in the extended space as follows
\cite{blanes10sac}
%
\begin{equation}\label{}
  \left\{ \begin{array}{l}
 u'=  A(t_1, u) \\
 t_1' = 1\\
 t_2' = 0
  \end{array} \right.
 \qquad \mbox{ and } \qquad
  \left\{ \begin{array}{l}
 u'=  B(t_2, u)\\
 t_1' = 0\\
 t_2' = 1.
  \end{array} \right.
\end{equation}

A more convenient way to split the system which transforms the
non-autonomous problems into autonomous is the following
\begin{equation}\label{eq.StandardSplitNA}
  \left\{ \begin{array}{l}
 u'=  A(t_1, u) \\
 t_1' = 0\\
 t_2' = 1
  \end{array} \right.
 \qquad \mbox{ and } \qquad
  \left\{ \begin{array}{l}
 u'=  B(t_2,u)\\
 t_1' = 1\\
 t_2' = 0.
  \end{array} \right.
\end{equation}
Notice that the explicit time-dependency in $A$ and $B$ is frozen
in each subproblem and the formal solution corresponds to the
exponential of the Lie operators where the time dependency in the
operators are frozen on each time interval. Unfortunately, to use
splitting methods with complex coefficients for non-autonomous
problems requires, in general, to compute $A(t,u),B(t,u)$ for
$t\in\mathbb{C}$, leading, in general, to badly conditioned
algorithms.

In this work we show that splitting method having one set of
coefficients real and positive valued, i.e.
\[
  a_i\in \mathbb{R}^+, \qquad b_i\in\mathbb{C}^+, \qquad
 (\mbox{or} \quad  a_i\in \mathbb{C}^+, \qquad
 b_i\in\mathbb{R}^+),
\]
allow to build algorithms where the operators $A(t,u),B(t,u)$ are
evaluated only for $t\in\mathbb{R}$, leading to well defined
methods. Several splitting methods with this structure have
already been constructed\footnote{ In \cite{castella09smw}, a
fourth-order method was obtained with $a_i\in\mathbb{R}^{+}$. In a
similar way, in \cite{blanes13oho} sixth-order schemes were also
explored with $a_i\in\mathbb{R}^{+}$. The coefficients can be
found at:
\texttt{http://www.gicas.uji.es/Research/splitting-complex.html}.
}.


We will also explore the case in which $\|B\|\ll\|A\|$, which we
refer as a perturbed problem. We first show how this class of
methods has to be used in these problems and next we study how to
build high order efficient methods for these problems.

\section{Splitting methods for non-autonomous problems}


Suppose we have a splitting method with say, $a_i\in\mathbb{R}^+$
and $b_i\in\mathbb{C}^+$. To solve the eq. (\ref{eq.Split_tDep})
we propose to take the time as one new coordinate and split the
system as follows \cite{blanes10sac}
\begin{equation}\label{eq:split1}
  \left\{ \begin{array}{l}
 u'=  A(t_1, u) \\
 t_1' = 1
  \end{array} \right.
 \qquad \mbox{ and } \qquad
  \left\{ \begin{array}{l}
 u'=  B(t_1, u)\\
 t_1' = 0.
  \end{array} \right.
\end{equation}
Let us denote by $\Phi_A^{[a_ih]}$ the map associated to the exact
solution (or a sufficiently accurate numerical approximation) of
the non-autonomous equation
\[
   \frac{d u}{dt}= {   } A(t, u), \qquad \quad
   t\in[t_n+c_{i-1}h,t_n+c_{i}h]
\]
with
\[
  c_i=\sum_{j=0}^{i} a_j,
\]
and $a_0=0$. Then, the splitting method (\ref{eq:splittingmethod})
for the non-autonomous equation reads now

\begin{eqnarray} \label{eq:splittingNoAut}
 \Psi(h) = \e^{h b_{m+1} L_{B_m}} \, \Phi_A^{[a_mh]} \,\cdots \,
  \e^{h b_2 L_{B_1}} \,
 \Phi_A^{[a_1h]} \, \e^{h b_1 L_{B_0}},
\end{eqnarray}
where
$
  B_i=B(t_n+c_ih,u).
$ 
 Notice that in this scheme, since $t_1$ is advanced with the
coefficients $a_i$ (and then it takes real values) the operators
$A(t,u)$ and $B(t,u)$ are evaluated on real values of $t$. On the
other hand, if $A(t,u)$ is an unbounded operator and a numerical
methods is used to approximate the flow $\Phi_A^{[h]}$, it must be
well defined for $0\leq h<h^*$ for some positive $h^*$, and this
is not guaranteed for general methods. For example, some
commutator-free methods up to fourth-order can be used. Given the
equation
\[
   \frac{d u}{dt}= A(t, u), \qquad \qquad  t\in[0,h]
\]
we have that
\[
  \Phi_A^{[h]} =  \e^{h L_{A(t_n+h/2,u)}}
\]
corresponds to a symmetric second order method, and
\begin{eqnarray}
  \Phi_A^{[h]} & = &  \e^{\frac{h}{2} (\alpha L_{A_1}+\beta L_{A_2})} \,
    \e^{\frac{h}{2} (\beta L_{A_1(u)}+\alpha L_{A_2(u)})}
    \nonumber \\
     & = & \Phi_{\frac{1}{2} (\beta A_1(u)+\alpha A_2(u))}^{[h]} \circ \Phi_{\frac{1}{2} (\alpha A_1(u)+\beta A_2(u))}^{[h]} \label{eq.CF4}
\end{eqnarray}
with $\alpha=\frac12-\frac{\sqrt{3}}{3}, \ \beta=1-\alpha$ and
$A_1(u)=A\left(t_n+\left(\frac12-\frac{\sqrt{3}}{6}\right)h,u\right), \
A_2(u)=A\left(t_n+\left(\frac12+\frac{\sqrt{3}}{6}\right)h,u\right)$,
%
corresponds to a fourth-order method
\cite{blanes06moan,thalhammer06afo}.
Notice that the Lie operators, since being derivatives, are
written in the  reverse order than the maps, and this is very
important to keep inmind for non-linear non-autonomous problems in
order to apply the method correctly (see \cite{blanes01smf} for
more details on the Magnus series expansion and Magnus integrators
for non-autonomous non-linear differential equations). This scheme
corresponds to the composition of the 1-flow maps for the
equations
\begin{eqnarray}
  u_1' & = &  \frac{h}{2} (\alpha A_1(u_1)+\beta A_2(u_1)), \qquad
    u_1(0)=u_0
    \nonumber \\
  u_2' & = &  \frac{h}{2} (\beta A_1(u_2)+\alpha A_2(u_2)), \qquad
    u_2(0)=u_1(1).
\end{eqnarray}
and the solution given by the map corresponds to $u_2(1)$.
 These second and fourth-order
commutator-free methods can be used for unbounded operators and
higher order commutator-free Magnus integrators for unbounded
operators are under investigation at this moment
\cite{blanes13wip}.

\subsection{Splitting methods for non-autonomous perturbed systems}

In some cases, the system can be considered as the perturbation of
an exactly solvable problem. In those cases, it is usually
convenient to split into the dominant part and the perturbation
and to build methods which take advantage of this relevant
property. However, if the problem is non-autonomous and the
time-dependency is not treated properly, the performance of the
methods designed for perturbed problems deteriorate considerably.

Suppose that $\|B(t,u)\|\ll\|A(t,u)\|$. To make this fact
more evident, we
replace $B$ by $\varepsilon B$
with
$|\varepsilon|\ll 1$\footnote{In most cases, this split is also
convenient for not necessarily very small perturbations, say
$\varepsilon<1/2$.}. In the autonomous case, for example, the
Lie-Trotter composition for this split satisfies
\begin{equation} \label{eq:strang2}
 \e^{h( L_{A+\varepsilon B})} = \e^{h\, L_A}
\e^{h \varepsilon L_{B}} +
 \frac12 \varepsilon h^2 [L_A,L_{B}] + \mathcal{O}(\varepsilon h^3),
\end{equation}
i.e., it has a local error of order $\mathcal{O}(\varepsilon
h^2)$.

Since $A$ and $B$ are qualitatively different for perturbed
problems, it is usual to consider $ABA$ and $BAB$ compositions.

An $m$-stage symmetric $BAB$ compositions given by
\begin{eqnarray} \label{eq:BAB}
\Psi(h) = \e^{h b_{m+1}\varepsilon L_B} \, \e^{h a_m L_A} \,\cdots \,
\e^{h b_2\varepsilon L_B} \, \e^{h a_1 L_A} \, \e^{h b_1\varepsilon
L_B},
\end{eqnarray}
with $a_{m+1-i}=a_i, \ b_{m+2-1}=b_i, \ i=1,2,\ldots$,  and $ABA$
compositions are given by
\begin{eqnarray} \label{eq:ABA}
\Psi(h) = \e^{h a_{m+1} L_A} \, \e^{h b_m \varepsilon L_B} \,\cdots \,
\e^{h a_2 L_A} \, \e^{h b_1\varepsilon L_B} \, \e^{h a_1 L_A},
\end{eqnarray}
with $a_{m+2-i}=a_i, \ b_{m+1-1}=b_i, \ i=1,2,\ldots$. We will use
the following short notation for these methods
\[
  (b_{m+1}, a_m , \,\cdots \, , b_2, a_1, b_1)
  \qquad \mbox{and} \qquad
  (a_{m+1}, b_m , \,\cdots \, , a_2, b_1, a_1).
\]

Notice that eq. (\ref{eq:BAB}) is a $BAB$ composition which, for
the particular case where $b_1=b_{m+1}=0$ transforms into a $ABA$
composition (but with a different computational cost). It seems
then natural to consider separately the following four cases:
\begin{enumerate}
  \item \label{BABaRbC} $BAB: \quad a_i\in \mathbb{R}^+, \quad
  b_i\in\mathbb{C}^+$,
  \item \label{BABaCbR} $BAB: \quad a_i\in \mathbb{C}^+, \quad
  b_i\in\mathbb{R}^+$,
  \item \label{ABAaRbC} $ABA: \quad a_i\in \mathbb{R}^+, \quad
  b_i\in\mathbb{C}^+$,
  \item \label{ABAaCbR} $ABA: \quad a_i\in \mathbb{C}^+, \quad  b_i\in\mathbb{R}^+$.
\end{enumerate}


The cases \ref{BABaCbR} and \ref{ABAaCbR} require to split the
system as follows
\begin{equation}\label{eq:splitPert1}
  \left\{ \begin{array}{l}
 u'=  A(t_1, u) \\
 t_1' = 0
  \end{array} \right.
 \qquad \mbox{ and } \qquad
  \left\{ \begin{array}{l}
 u'= \varepsilon  B(t_1, u)\\
 t_1' = 1,
  \end{array} \right.
\end{equation}
so $t_1$ will take real values. This split is similar to the one
shown in the previous section by changing the roles of $A$ and
$B$.

In the extended phase space these two systems are equivalent to
solve separately the following system written in terms of Lie
operators
\begin{equation}\label{eq:splitPert2}
  \frac{d}{dt} \left\{ \begin{array}{l}  u \\ t_1
  \end{array} \right\} =
   \underbrace{\left(  {   } A(t_1, u) \frac{\partial }{\partial u}
  + 0 \cdot\frac{\partial }{\partial t_1} \right)}_{\mathcal{A}}
  \left\{ \begin{array}{l}
 u \\ t_1
  \end{array} \right\}
\end{equation}

\begin{equation}\label{eq:splitPert2}
  \frac{d}{dt} \left\{ \begin{array}{l}  u \\ t_1
  \end{array} \right\} =
  \underbrace{\left( \varepsilon {   } B(t_1, u) \frac{\partial }{\partial u}
  + 1\cdot\frac{\partial }{\partial t_1} \right)}_{\mathcal{B}}
  \left\{ \begin{array}{l}
 u \\ t_1
  \end{array} \right\}.
\end{equation}
The commutators of the Lie operators ${\mathcal{A}}$ and
${\mathcal{B}}$, which measure the error of the splitting methods
in the extended phase space is
\[
 [h{\mathcal{A}},h{\mathcal{B}}]=
 h^2\big({\mathcal{A}}{\mathcal{B}}-{\mathcal{B}}{\mathcal{A}}\big)=
  h^2\left( \varepsilon [{   } A, {   } B] - \frac{d{   } A(t_1,u)}{d t_1} \right)
  \frac{\partial }{\partial u} = \mathcal{O}(h^2)
\]
which is not proportional to $\varepsilon$ due to the term
$\frac{d A(t_1,u)}{d t_1}$, and this also happens with higher
order commutators.

The cases \ref{BABaRbC} and \ref{ABAaRbC} are associated to the
split
\begin{equation}\label{eq:splitPert3}
  \left\{ \begin{array}{l}
  u'=  A(t_1, u) \\
 t_1' = 1
  \end{array} \right.
 \qquad \mbox{ and } \qquad
  \left\{ \begin{array}{l}
  u'= \varepsilon  B(t_1, u)\\
 t_1' = 0.
  \end{array} \right.
\end{equation}
This system can be written in the extended phase space as
\begin{equation}\label{eq:splitPert4}
  \frac{d}{dt} \left\{ \begin{array}{l}  u \\ t_1
  \end{array} \right\} =
   \underbrace{\left(  {   } A(t_1, u) \frac{\partial }{\partial u}
  + 1 \cdot\frac{\partial }{\partial t_1} \right)}_{\mathcal{A}}
  \left\{ \begin{array}{l}
 u \\ t_1
  \end{array} \right\}
\end{equation}

\begin{equation}\label{eq:splitPert5}
  \frac{d}{dt} \left\{ \begin{array}{l}  u \\ t_1
  \end{array} \right\} =
  \underbrace{\left( \varepsilon {   } B(t_1, u) \frac{\partial }{\partial u}
  + 0\cdot\frac{\partial }{\partial t_1} \right)}_{\mathcal{B}}
  \left\{ \begin{array}{l}
 u \\ t_1
  \end{array} \right\}
\end{equation}
where now
\[
 [h{\mathcal{A}},h{\mathcal{B}}]=
  \varepsilon h^2\left( [{   } A, {   } B] + \frac{d B(t_1,u)}{d t_1} \right)
  \frac{\partial }{\partial u} = \mathcal{O}(\varepsilon h^2)
\]
which is proportional to the small parameter $\varepsilon$ (see
\cite{blanes10sac} for more details).

Obviously, this split makes sense if one can exactly solve the
non-autonomous equation associated to the dominant part
\begin{equation}\label{}
 \frac{d u}{dt}=  A(t,u)
\end{equation}
at a relatively low computational cost (or one can numerically
solve it up to sufficiently high accuracy and at a relatively low
computational cost) being the commutator-free Magnus integrators
an appropriate choice in most cases. A similar methods was used in
\cite{bader11fmf} for perturbed Schr\"odinger and Gross-Pitaevskii
equations, but in those problems negative real coefficients are allowed.

\subsection{Order conditions}

 For consistent symmetric methods
we can formally write
\begin{eqnarray}
 \Psi(h) &=&
 \exp\left(h(L_A+\varepsilon L_B)+h^3 \, \left( \eps \, p_{aba} [[L_A,L_B],L_A] +
 \eps^2\, p_{abb} [[L_A,L_B],L_B]) \right)  \right. \nonumber \\
 & & \qquad \left. + h^5 \left( \eps\, p_{abaaa} [[[[L_A,L_B],L_A],L_A],L_A]  + \mathcal{O}(\eps^2) \right)+
\mathcal{O}(\varepsilon h^7)\right). \label{eq:ordereqtsn}
\end{eqnarray}

Following \cite{blanes13nfo}, for the composition (\ref{eq:ABA})
we have that (taking $m=s$)
\begin{eqnarray}
p_{aba} & \sim & \frac12 \sum_{i=1}^{s} b_{i} c_{i} (1-c_{i})-\frac{1}{12}
,\qquad p_{abb}  \sim   \sum_{i=1}^{s} \frac{1}{2} b_{i}^{2} \, c_{i} +
\sum_{1 \leq i < j \leq s} b_{i} b_{j} c_{j}-\frac{1}{3},  \nonumber \\
p_{abaaa} & \sim &  \sum_{i=1}^{s} b_i \, c_i^{4}-\frac{1}{5},
 \label{p_abaaa}
\end{eqnarray}
%
with $a_0=0$ and $c_{s+1}=1$.
The symbol $\sim$ indicates that, if the low order conditions are satisfied,
both terms are proportional so, if the r.h.s of $p_{aba}$ and
$p_{abb}$ vanish then $p_{aba}=p_{abb}=0$.
%
 Here, the polynomial $p_{abaaa}$
corresponds to the dominant error term in fourth-order methods for
perturbed problem. This algebraic analysis remains also valid for unbounded operators under appropriate conditions on the operators (see \cite{hansen08esf} for more details).

If we take $a_1=a_{m+1}=0$ we obtain a $BAB$ composition, and the
equations  of $p_{aba}$ and $p_{abb}$ can be easily adjusted to
obtain $BAB$ compositions.



Second order symmetric methods which cancel the terms of order
$h^{2p+1}\varepsilon$ for $p=1,2,\ldots,s$ and for different
values of $s$ exist with positive and real coefficients
\cite{mclachlan95cmi}. The error of these methods is of order
$\mathcal{O}(h^{2s+1}\varepsilon + h^3\varepsilon^2)$ and we say
the methods have effective order $(2s,2)$. For instance, a method
which satisfies $p_{aba}=p_{abaaa}=0$ has effective order (6,2),
and this can be attained with the sequence \cite{mclachlan95cmi}
\begin{equation}\label{eq.62}
  (b_{1}, a_1  , b_2, a_2 , b_2, a_1,
  b_1)=\left(\frac1{12},\frac{5-\sqrt{5}}{10},\frac5{12},\frac1{\sqrt{5}},
  \frac5{12},\frac{5-\sqrt{5}}{10},\frac1{12}\right).
\end{equation}

Fourth-order methods require to satisfy
$ 
 p_{aba}=0, \ \  p_{abb}=0,
$ 
and this can not be accomplished with $a_i,b_i$ real and positive
valued coefficients. We are then interested on the existence of
methods in which $a_i\in\mathbb{R}^+$ and $b_i\in\mathbb{C}^+$. To
get splitting methods where the coefficients satisfy these
constraints we fix the values of the coefficients $a_i\in(0,1)$
such that consistency and symmetry is satisfied, and leave the
coefficients $b_i$ to solve the order conditions.
Obviously, since the coefficients are chosen real and positive,
the equations only admit complex solutions for the coefficients
$b_i$. Among all solutions obtained we will choose solutions with
positive real part, i.e. $b_i\in\mathbb{C}^+$ from the set of all
solutions found (in case these solutions exist).


Let us now analyse the number of free parameters and computational
cost of $ABA$ and $BAB$ compositions in order to choose the most
appropriate sequence: A symmetric $(2k)$-stage $BAB$ composition
has $k$ coefficients $a_i$ and $k+1$ coefficients $b_i$ while an
$ABA$ sequence has $k+1$ coefficients $a_i$ and $k$ coefficients
$b_i$ so, the $BAB$ composition has one more free parameter to
solve the equations.
%
In addition, since the dominant part is associated to the
coefficients $a_i$
and requires the numerical solution of a non-autonomous
differential equation,
%
it is not usually possible to concatenate the last map
in one step with the first one in the following step, and in
practice an $ABA$ composition with the same number of stages
as a $BAB$ composition can be computationally more costly up to
one additional stage. For these reasons (number of free parameters to solve the equations and the computational cost) we only consider in this work
$BAB$ compositions.

\subsection{Fourth-order methods}

Fourth-order methods can be obtained with a 4-stage composition
\begin{equation}\label{eq:composBAB4}
   (b_1 \, a_1 \, b_2 \, a_2 \, b_3 \, a_2 \, b_2 \, a_1 \, b_1)
\end{equation}
which satisfy the consistency conditions $a_1+a_2=1/2, \ \
2(b_1+b_2)+b_3=1. $
We can fix the values of $a_1$ such that $a_1\in(0,1/2)$ and take,
e.g. $b_1,b_2$ to solve the equations $p_{aba}=p_{abb}=0$.
The choice
$a_1=\frac14$ leads to the solution obtained in
\cite{castella09smw}. However, we can take $a_1$ as a free
parameter to minimise the dominant error term\footnote{We minimise
the real part of the dominant error because after each time step
we will remove the imaginary part of the numerical solution, i.e.
$u_{n+1}=Re(\Psi(h)u_n)$.}
\begin{equation}\label{eq:Min_pabaaa}
  \min_{a_1\in(0,1/2)} \left| Re(p_{abaaa})\right| =
  \min_{a_1\in(0,1/2)} \left|\sum_{i=1}^{4} Re(b_i) \,
  c_i^{4}-\frac{1}{5}\right|.
\end{equation}

The corresponding system of polynomial equations with two unknowns
have only two solutions (complex conjugate to each other) for each
choice of $a_1$ and with this process we obtain following method:
    \begin{equation}
    \begin{aligned}\label{eq:newcoeff}
     b_1 =\ &0.018329102861074364-0.10677008344599524i,\\
     a_1 =\ & 0.13505265889288437,\\
     b_2 =\ & 0.2784394345454581+0.20041452008768607i,\\
     a_2 =\ & 0.36494734110711563,\\
     b_3 =\ & 0.40646292518693505-0.18728887328338165i.\\
        \end{aligned}
        \end{equation}

A 5-stage $BAB$ composition 
has the same number of coefficient $b_i$ and for this reason
we have not considered it. To vanish the dominant error term at
order 6 we need at least a 6-stage composition
\begin{equation}\label{eq:composBAB6}
   (b_1 \, a_1 \, b_2 \, a_2 \, b_3 \, a_3 \, b_4 \, a_3 \, b_3 \, a_2 \, b_2 \, a_1 \,
   b_1)
\end{equation}
where the coefficients $b_i$ are used to satisfy the conditions
$
 p_{aba}=p_{abb}=p_{abaaa}=0.
$ 
in addition to consistency.

The goal of this work is not to make an exhaustive search of
methods but to show this class of methods are of interest for
non-autonomous problems and to indicate how highly efficient
methods could be obtained, and the optimal method can depend on
the algebraic structure of each problem\footnote{Higher order and
more efficient methods require a considerably deeper analysis, and
methods belonging to this class as well as more general methods
are being considered by the authors of Ref. \cite{blanes13oho}.}.
Then, just as an illustration we take $ a_1=a_2=a_3=\frac16$.
We have obtained one complex solution (and its complex conjugate)
with coefficients:
    \begin{equation}
    \begin{aligned}\label{eq:newcoeff1}
     a_1 =\ &  a_2 = a_3 = 1/6,\\
     b_1 =\ & 0.05753968253968254 - 0.007886748775536424i,\\
     b_2 =\ & 0.20476190476190473 + 0.04732049265321855i,\\
     b_3 =\ & 0.16309523809523818 - 0.11830123163304637i,\\
     b_4 =\ & 0.14920634920634912 + 0.15773497551072851i.\\
        \end{aligned}
        \end{equation}

\section{Numerical examples} \label{sect:nt}

To analyse the performance of the new methods we first consider a
simple non-autonomous ODE as a test bench of the methods and next
we apply the methods to a linear non-autonomous PDE and a
non-linear non-autonomous PDE. We compare the performance of the
methods with complex coefficients versus other methods which
involve real coefficients. We choose
the
(6,2) splitting method (\ref{eq.62}) which is a method of second
order. As a fourth-order method we consider extrapolation (which
involves substraction of quantities) where the Strang splitting
symmetric second order method is used as the basic scheme to raise
the order. To be more precise, we consider
\begin{equation} \label{eq:strang2b}
 S(h)= \e^{h/2\, L_{B_1}} \, \e^{h L_{A_0}} \, \e^{h/2\, L_{B_0}}
\end{equation}
where we denote $A_i=A(t_n+ih,u), \ B_i=B(t_n+ih,u)$. This scheme can
be considered as the standard  Strang decomposition applied to the
non-autonomous system, but if we split it as shown in
(\ref{eq.StandardSplitNA}). If we take $S(h)$ as the basic method,
high order methods by extrapolation can be obtained and they only
involve positive time steps. A fourth-order method is given by the
composition
\begin{equation}\label{eq:extap4a}
  \Phi^{[4]}(h) = \frac43 S\left(\frac{h}{2}\right) S\left(\frac{h}{2}\right) - \frac13 S(h)
\end{equation}
which in our case it can be written as
\begin{equation}\label{eq:extap4a2}
  \Phi^{[4]}(h) =
  \frac43  \e^{\frac{h}4\, L_{B_1}} \, \e^{\frac{h}2 L_{A_{1/2}}}
  \, \e^{\frac{h}2\, L_{B_{1/2}}}   \,
  \e^{\frac{h}2 L_{A_{0}}} \, \e^{\frac{h}4\, L_{B_0}} -
   \frac13 \e^{\frac{h}2\, L_{B_{1}}} \, \e^{h L_{A_0}} \, \e^{\frac{h}2\, L_{B_0}}.
\end{equation}

The following schemes with real coefficients are then considered:
\begin{itemize}
\item {\bf Strang}: The second-order symmetric Strang splitting method (as a reference method);
\item {\bf (6,2)}: The symmetric splitting method of effective order (6,2) whose
   coefficients are given in (\ref{eq.62});
\item {\bf (EXT4)}: The fourth-order extrapolation method (\ref{eq:extap4a2});
\end{itemize}
and the following schemes with real coefficients and
$a_i\in\mathbb{R}^+$ are considered:
\begin{itemize}
\item {\bf (RC4)}: The 4-stage fourth-order method from \cite{castella09smw};
\item {\bf (O4)}: The 4-stage fourth-order method built in \cite{blanes13oho},
 whose coefficients are available at
 \texttt{http://www.gicas.uji.es/Research/splitting-complex.html}, and referred as
 "Order 4 (optimized)";
\item {\bf (SM4)}: The new optimized 4-stage fourth-order method given in (\ref{eq:newcoeff});
\item {\bf (SM(6,4))}: The new 6-stage fourth-order method whose coefficients are given in (\ref{eq:newcoeff1});
\end{itemize}

The numerical approximations  $u_n$ obtained by a given method,
$\Psi(h)$, which involve complex coefficients are computed as $u_n
= \Re(\Psi(h) u_{n-1})$, i.e. we project on the real axis after
completing each time step. To measure the performance of the
methods we compute the error of each method at the end of the time
integration (we take as the exact solution a numerical
approximation computed to a high precision) and we take as the
cost of the method the number of evaluations of $\Phi^{[h]}_A$
which usually carries most of the computational cost.

\paragraph*{\bf Example 1} Let us consider the non-autonomous and
non-linear perturbed equation
\begin{equation}
  q'' + \Omega(t)^2 q = -\varepsilon
  \sum_{j=1}^s\sin (q-\omega_j t), \qquad q\in \mathbb{R}.  \label{wave1}
\end{equation}
When $\Omega$ is a constant, the system
describes the motion of a charged particle in a magnetic field
perturbed by $s$ electrostatic plane waves, each with the same
wavenumber and amplitude, but with different temporal frequencies
$\omega_{j}$ \cite{candy91asi}. This equation can be written as a
first order system of equations
\[
  \frac{d}{dt}
   \left\{ \begin{array}{c} q \\ p \end{array} \right\} =
   \left( \begin{array}{cc} 0 & 1 \\ -\Omega(t)^2 & 0 \end{array} \right)
   \left\{ \begin{array}{c} q \\ p \end{array} \right\} +
   \varepsilon \left\{ \begin{array}{c} 0 \\ -
  \sum_{j=1}^s\sin (q-\omega_j t) \end{array} \right\}
\]
which we split as follows
\[
  \frac{d}{dt}
   \left\{ \begin{array}{c} q \\ p \end{array} \right\} =
   \left( \begin{array}{cc} 0 & 1 \\ -\Omega(t_1)^2 & 0 \end{array} \right)
   \left\{ \begin{array}{c} q \\ p \end{array} \right\},  \qquad
   \quad
  \frac{dt_1}{dt} = 1
\]
and
\[
  \frac{d}{dt}
   \left\{ \begin{array}{c} q \\ p \end{array} \right\} =
   \varepsilon \left\{ \begin{array}{c} 0 \\ -
  \sum_{j=1}^s\sin (q-\omega_j t_1) \end{array} \right\}.
\]
The linear part has, in general, no solution in closed form and we
approximate its flow using the the 4th-order commutator-free Magnus
integrator (\ref{eq.CF4}) which for this problem reads\footnote{Here, the method is written in terms of exponentials of matrices, i.e. maps, so they appear in the reverse order as the Lie operators in (\ref{eq.CF4}).  }
\begin{eqnarray}\label{eq:Rotation}
 \Phi^{[a_i h]}_A & = &
  \e^{\frac{a_ih}{2} (\beta A_1+ \alpha A_2)}
 \,
  \e^{\frac{a_i h}{2} (\alpha A_1+ \beta A_2)}  \\
 & = &
  \exp\left[ \frac{a_i h}{2} \left(
   \begin{array}{cc}
     0 & 1  \\
     -(\beta \Omega_1^2+ \alpha \Omega_2^2) & 0
   \end{array}
  \right) \right]
    \
  \exp\left[ \frac{a_i h}{2} \left(
   \begin{array}{cc}
     0 & 1  \\
     -(\alpha \Omega_1^2+\beta  \Omega_2^2) & 0
   \end{array}
  \right) \right]
    \nonumber
\end{eqnarray}
where $\Omega_i=\Omega(t_n+c_ih)$ and the exponential of each matrix can be easily computed taking
into account that
\begin{equation}\label{eq:Rotation}
  \exp\left[ \tau \left(
   \begin{array}{cc}
     0 & 1  \\
     -\Omega^2 & 0
   \end{array}
  \right) \right]=  \left(
   \begin{array}{cc}
     \cos(\tau\Omega) & \frac1{\Omega}\sin(\tau\Omega)  \\
     -\Omega\sin(\tau\Omega) & \cos(\tau\Omega)
   \end{array}
  \right).
\end{equation}
The evolution for the perturbation is immediate since both $q$ and
$t$ are frozen.

Notice that if the computational cost is dominated by the
evaluation of the time-dependent functions and the rotation
matrix, then since $a_i\in \mathbb{R}$ the overall cost does not
change considerably either if $b_i$ is real or complex.

\begin{figure}[!t]
\makebox{\psfig{figure=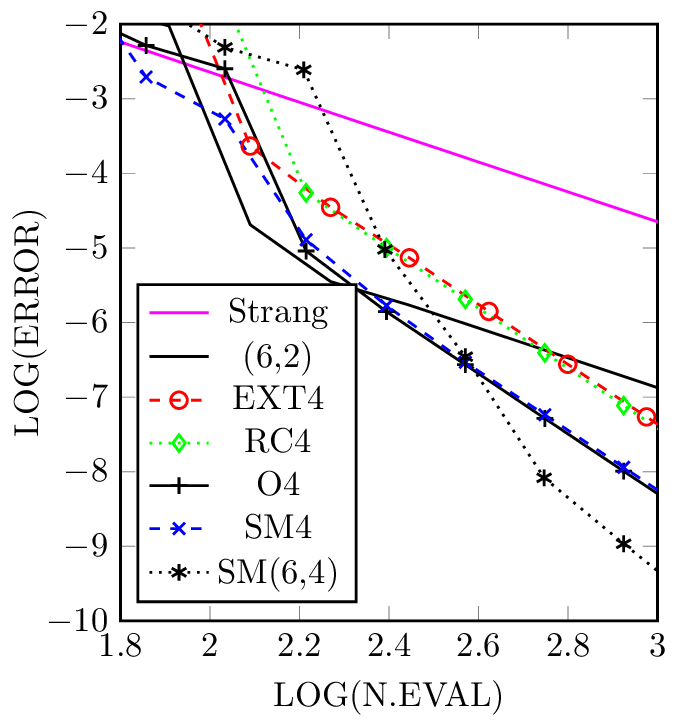,height=7.2cm,width=7cm}}
\makebox{\psfig{figure=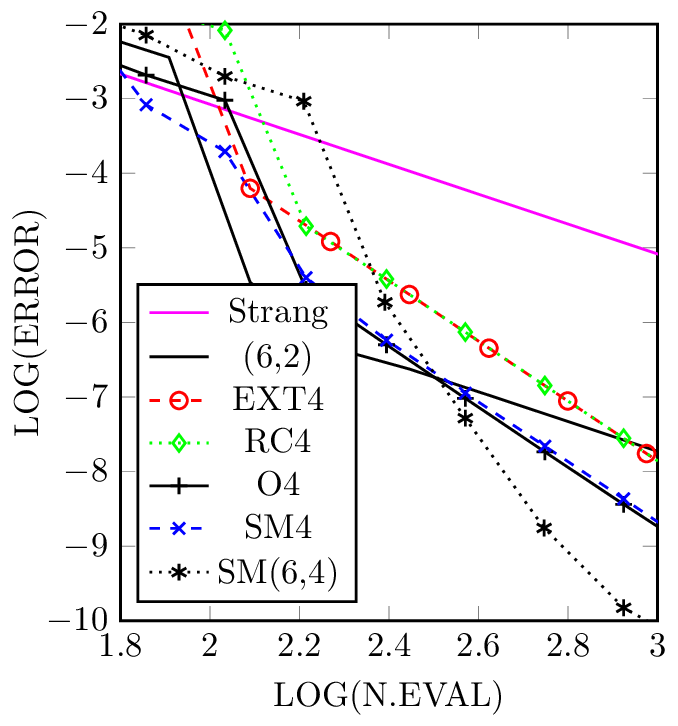,height=7.2cm,width=7cm}}
 \caption{ Error versus number of evaluations of $\Phi^{[h]}_A$ for the numerical
    integration in Example 1 at $t=2\pi$ for $\varepsilon=\frac14$
     (left panel) and $\varepsilon=\frac1{10}$ (right panel).}
\label{fig:phosc}
\end{figure}

For the numerical experiments we take
$\Omega(t)=1+\frac{1}{2}\cos(\frac{3}{2}t)$ and the same initial
conditions and parameters as given in \cite{candy91asi}:
$
 (q_0,p_0,t_0)=(0,11.2075,0),  \ s=3
 \ \omega_j=7j.
$ 
We integrate for $t\in[0,2\pi]$ and measure the error at the final
time. All the computations are done for $\varepsilon=1/4$ and
$\varepsilon=1/10$. Fig.~\ref{fig:phosc} shows the error versus
the number of evaluations for different methods. We clearly
observe the superiority of the methods which consider complex
coefficients versus the lower order splitting methods with real
coefficients or extrapolation when high accuracy is desired as
well as the high performance of the new methods.

\begin{figure}[!t]
\begin{center}
\makebox{\psfig{figure=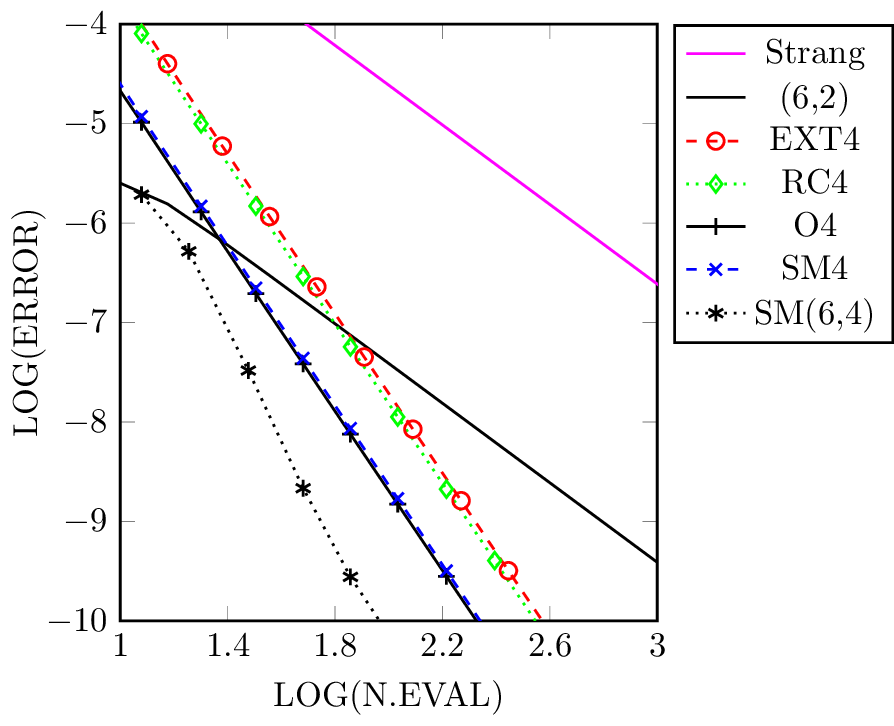,height=7.2cm,width=10cm}}
\end{center}
 \caption{  Efficiency comparison between algorithms for the linear parabolic
    equation (\ref{eq:lrd}) with parameters $\mu=1/6, w=2$ at final time $t=1$.}
\label{fig:llnppde}
\end{figure}

\paragraph*{\bf Example 2: A linear parabolic equation.}
%
The next test-problem is the following scalar parabolic equation
in one-dimension
\begin{eqnarray} \label{eq:lrd}
\frac{\partial u(x,t)}{\partial t} = \alpha(t)^{2} \Delta u(x,t) +
V(x,t) u(x,t),  \qquad u(x,0) = u_{0}(x),
\end{eqnarray}
with $u_{0}(x) = \sin(2 \pi x)$ and periodic boundary conditions
in the space domain $[0,1]$. We take $\alpha(t) = \frac{1}{4}+\mu
\cos(w t)$,
 $V(x,t)=\frac1{10}\left(3(1-e^{-t})+\sin(2 \pi x)\right)$
and discretize in space
$$x_j= j (\delta x), \qquad j=1,\ldots,N \quad \mbox{ with } \quad \delta x = 1/N,
$$
thus arriving at the differential equation
\begin{equation} \label{eq:problem0}
\frac{dU}{dt} = \alpha(t)^{2} A U + B(t) U,
\end{equation}
where $U=(U_1, \ldots,U_N)=(u_1, \ldots,u_N) \in \R^N$. The
Laplacian $\Delta$ has been approximated by the matrix $A$ of size
$N\times N$ given by\footnote{Our main purpose here is just to
illustrate the performance of the new splitting methods. In this
sense, the particular scheme used to discretize in space is
irrelevant.  For that reason, and to keep the treatment as simple
as possible, we have applied a simple second-order finite
difference scheme in space.}
\begin{equation}\label{eq.MatrixA}
A=\frac{1}{(\delta x)^2}\begin{pmatrix}
-2&1& & & 1\\
1 & -2 & 1 \\
 & 1 & -2 & 1 \\
 & & \ddots & \ddots & \ddots \\
1 & & & 1 & -2\end{pmatrix},
\end{equation}
and  $B(t)=\mbox{diag}(V(x_1,t), \ldots, V(x_N,t) )$. We take
$\mu=1/6, \ w=2$, $N=100$ points and compare different
composition methods by computing the corresponding approximate
solution on the time interval $[0,1]$.
 We compute the $2$-norm error of the numerical solution with respect to the exact solution of the semidiscretised equation (computed numerically up to a sufficiently high accuracy) at time
$t=1$. The results are collected in
Fig.~\ref{fig:llnppde} where the superiority of the splitting
methods with complex coefficients is also manifest.

\paragraph*{\bf Example 3: The semi-linear reaction-diffusion equation of Fisher.}
%
Our final test-problem is the following non-linear parabolic
scalar equation in one-dimension
\begin{eqnarray} \label{eq:nlrd}
\frac{\partial u}{\partial t}  = \alpha(t)^{2}\Delta u + F(u,t),
\qquad u(x,0) = u_{0}(x),
\end{eqnarray}
with periodic boundary conditions in the space domain $[0,1]$. We
take, in particular, the Fisher's potential
\[
 F(u)=\gamma(t)u(1-u), 
\]
with $\gamma(t)=(2-e^{-\beta t})/100$ and $\alpha(t) = \frac{1}{4}+\mu \cos(w t)$.

The splitting considered here corresponds to solving, on one hand,
the linear equation with $A$ given by (\ref{eq.MatrixA}) and on
the other hand, the nonlinear ordinary differential equation
$$
\frac{\partial u}{\partial t} = \gamma(t)u(1-u)
$$
with initial condition
$
u(x,0)=u_{0}(x).
$
After discretization in space,
we arrive at the differential equation
\begin{equation} \label{eq:problem1}
\frac{dU}{dt} = \alpha(t)^{2} A U + F(U,t),
\end{equation}
where $U=(U_1, \ldots,U_N)=(u_1, \ldots,u_N) \in \R^N$, A is a
circulant matrix of size $N\times N$ as in the the previous linear
case and $F(U,t)$ is now defined by
$$
F(U,t)=\gamma(t)\big(U_1(1-U_1), \ldots, U_N(1-U_N) \big).
$$
Here we consider splitting technique (\ref{eq:splitPert3}) as
\begin{equation}\label{eq:splitPerttt}
  \left\{ \begin{array}{l}
 U'=  \alpha(t_{1})^{2}A U  \\
 t_1' = 1
  \end{array} \right.
 \qquad \mbox{ and } \qquad
  \left\{ \begin{array}{l}
 U'= \gamma_{1}U(1-U)\\
 t_1' = 0.
  \end{array} \right.
\end{equation}
where $\gamma_{1}=\gamma(t_{1})$. Since $\gamma(t)$ is frozen at
real values of $t$, it must be considered as a constant, and the
scalar equations can be solved analytically


$$
u(x,h) = u_{0}(x)\frac{\e^{\gamma_{1}h}}{1+u_{0}(x)
(\e^{\gamma_{1}h}-1)} \ ,
$$
which is well defined for small complex time $h$.
We proceed in the same way as for the previous linear case,
starting with $u_{0}(x) = \sin(2 \pi x)$.

\begin{figure}[!t]
\begin{center}
\makebox{\psfig{figure=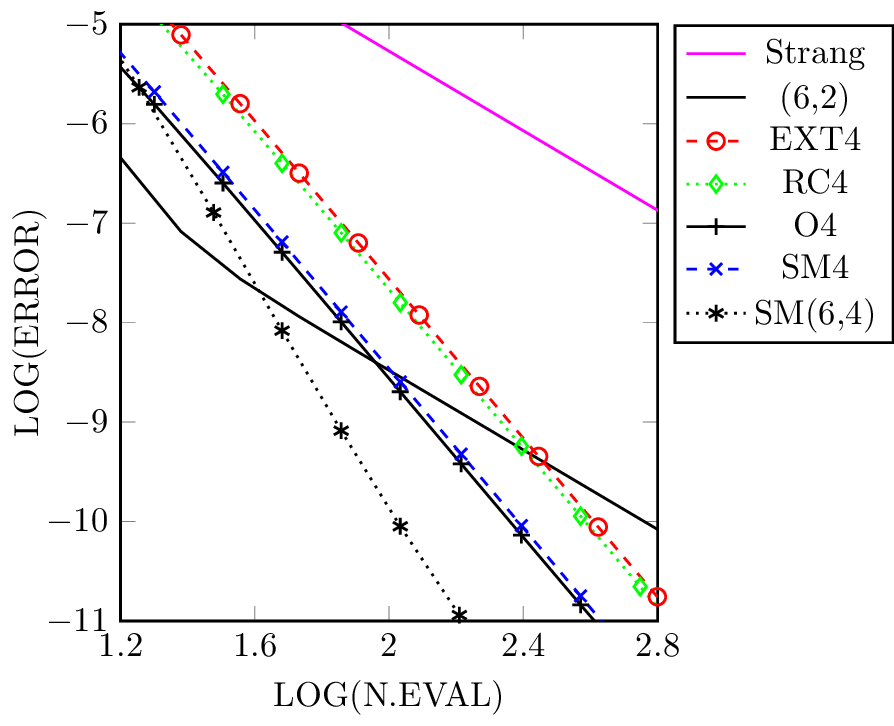,height=7.2cm,width=10cm}}
\end{center}
 \caption{  Efficiency comparison between algorithms for the equation of Fisher
    with parameters $\beta=1, \mu=1/6, w=2$ at final time $t=1$.}
\label{fig:fisssh}
\end{figure}

We choose $\beta=1,\mu=1/6, \ w=2$, $N=100$ and compute the error
at the
final time $t=1$ by applying the same composition methods as in
the linear case. The results are collected in Fig.~\ref{fig:fisssh}.

\section{Conclusions}

We have considered the numerical integration of non-autonomous
separable parabolic equations using high order splitting methods
with complex coefficients. A straightforward application of
splitting methods with complex coefficients to non-autonomous
problems require the evaluation of the time-dependent functions on
the operators at complex times, and the corresponding flows in the
numerical scheme are, in general, not well conditioned. To avoid
this trouble, in this work we consider a class of methods in which
one set of the coefficients belong to the class of real and
positive numbers. Taking the time as a new coordinate and an
appropriate splitting of the system allows us to build numerical
schemes where all time-dependent operators are evaluated at real
values of the time. This technique shows a great interest for
perturbed systems, and this problem is analysed in more detail. In
this case, the flow of the dominant part has to be advanced using
the real coefficients. We have analysed the algebraic structure of
the problem and the cost of the algorithms in order to build
efficient high order methods, and some few new methods are reported
as an illustration. Higher order and more efficient methods require a considerably deeper analysis, and methods belonging to this class as well as more general methods are being considered by the authors of Ref. \cite{blanes13oho} and will be published elsewhere.
Several numerical examples are considered where it is shown the
good performance of this class of methods. We have shown that
splitting methods  with complex coefficients can also be used on
non-autonomous non-linear parabolic problems and they can show a
good performance so, high order and more efficient schemes
following the guidelines presented in this work can be of great
interest. We must also remark that order reductions are expected
for problems with Dirichlet or Newman boundary conditions on
bounded domains, and the performance of high order methods on
these problems diminishes, being an interesting problem that needs
further investigation.

\section*{Acknowledgements}


The authors thank the referees for their suggestions to improve
the presentation of this work. The work of Sergio Blanes has been
supported by Ministerio de Ciencia e Innovación (Spain) under
project MTM2010-18246-C03 and the Ministerio de Educación, Cultura
y Deporte, under Programa Nacional de Movilidad de Recursos
Humanos del Plan Nacional de I-D+i 2008-2011 (PRX12/00547). The
work of Muaz Seydao\u{g}lu has been supported by the Turkish
Council of High Education through a grant to visit the Instituto
de Matem\'atica Multidisciplinar at University Polytechnique of
Valencia where this work was carried out.




\end{document}